\title{Hypermaps and multiply quasiplatonic Riemann surfaces}

\author{Gareth A. Jones\\
School of Mathematics\\
University of Southampton\\
Southampton SO17  1BJ, U.K.\\
{\tt G.A.Jones@maths.soton.ac.uk}
}

\documentclass[12pt]{article}
\usepackage[latin1]{inputenc}
\usepackage{a4wide}
\usepackage{latexsym}
\usepackage{amsmath}
\usepackage{amssymb}
\usepackage{amsfonts}
\newtheorem{thm}{Theorem}[section]
\newtheorem{lemma}[thm]{Lemma}
\newtheorem{cor}[thm]{Corollary}
\newtheorem{prop}[thm]{Proposition}
\date{}
\begin{document} 

\maketitle

\centerline{\bf This paper is dedicated to Antonio Mach\`\i, on his 70th birthday.}

\bigskip

\begin{abstract}
\noindent Generalising an example by Girondo and Wolfart, we use finite group theory to construct Riemann surfaces admitting two or more regular dessins (i.e.~orientably regular hypermaps) with automorphism groups of the same order, and in many cases with isomorphic automorphism groups.
\end{abstract}

% 05C10 = Planar graphs, geometric and topological aspects of graph theory
% 14H57 = Dessins d'enfants theory
% 20B25 = Finite automorphism groups of algebraic, geometric or combinatorial structures


\noindent{\bf Running head:} Hypermaps.

\section{Introduction}

Hypermaps were introduced by Cori~\cite{Cor} in 1975, and their theory was developed in the following years, mainly by Cori and Mach\`i~\cite{CoMa81, CoMa88, CoMa90, CMPV, Mach82, Mach84, Mach85}; see~\cite{CoMa} for a comprehensive survey. As originally defined, hypermaps are algebraic and combinatorial structures which model certain embeddings of finite graphs in compact oriented surfaces. Although the theory has subsequently been extended to include surfaces which may be non-orientable, with boundary~\cite{IS, JS93}, or non-compact~\cite{Jon08}, the main emphasis of current research still concerns hypermaps on compact oriented surfaces. Recently, much of this research has been motivated by the unexpected role such hypermaps play in Grothendieck's theory of {\em dessins d'enfants}~\cite{Gro}, where they provide a vital link between compact Riemann surfaces, algebraic curves, and the Galois theory of algebraic number fields; see~\cite{GG, JS93, JS96, Sin3, Wol2} for their connections with this theory. Another recent application (originally motivating this paper) is in the construction of Beauville surfaces~\cite[p.~159]{Bea}, rigid complex surfaces formed from pairs of regular hypermaps with isomorphic automorphism groups~\cite{Cat, FuJo}.

A hypermap $\cal H$, as originally defined in~\cite{Cor}, is a pair of permutations generating a finite transitive permutation group $G$. As such, it is equivalent to a permutation representation $\Delta\to G$ of a triangle group $\Delta$, acting on the cosets of a subgroup $K$ of finite index in $\Delta$. Now $\Delta$ acts as a cocompact group of automorphisms of a simply connected Riemann surface $\mathbb U$, and the quotient ${\mathbb U}/K$ is a compact Riemann surface $X$. By Bely\u\i's Theorem~\cite{Bel}, as reinterpreted by Grothendieck~\cite{Gro}, Wolfart~\cite{Wol2} and others, and as eventually proved by Wolfart~\cite{Wol1} and Koeck~\cite{Koe}, a compact Riemann surface is obtained in this way from a hypermap if and only if it is defined, as a projective algebraic curve, over the field $\overline{\mathbb Q}$ of algebraic numbers. Thus $\cal H$ carries with it the structure of a compact Riemann surface, or equivalently a projective algebraic curve, defined over $\overline{\mathbb Q}$; as such, we will call $\cal H$ a {\em dessin}, and $X$ a {\em Bely\u\i\/ surface}.

The most symmetric  Bely\u\i\/ surfaces are the {\em quasiplatonic surfaces}, those obtained from regular dessins, i.e.~orientably regular hypermaps $\cal H$; this is equivalent to $G$ being a regular permutation group, to $K$ being a normal subgroup of $\Delta$, and to the Bely\u\i\/ function ${\mathbb U}/K\to{\mathbb U}/\Delta\cong {\mathbb P}^1({\mathbb C})$ being a regular covering. As a dessin, $\cal H$ then has automorphism group ${\rm Aut}\,{\cal H}\cong G\cong\Delta/K$. (This is the orientation-preserving automorphism group of the hypermap, inducing conformal automorphisms of the Riemann surface $X$; as a hypermap, $\cal H$ may also have orientation-reversing automorphisms, acting anti-conformally on $X$.)

It is well known that a quasiplatonic Riemann surface $X$ may support two (or more) regular dessins ${\cal H}_1$ and ${\cal H}_2$; such surfaces are called {\em multiply quasiplatonic}. In the most familiar cases these dessins correspond to an inclusion $\Delta_1<\Delta_2$ between triangle groups, with a normal subgroup $K$ of $\Delta_1$ also normal in $\Delta_2$, inducing an inclusion $G_1<G_2$ between the corresponding automorphism groups. Inclusions between triangle groups have been classified by Singerman~\cite{Sin2}, and the corresponding relationships between dessins have been studied in various papers by Girondo, Torres and Wolfart, such as~\cite{Gir, GTW, GW}.

Instead, we will be concerned here with a less common situation, where a multiply quasiplatonic surface $X$ supports regular dessins ${\cal H}_i$ with automorphism groups $G_i$ of the same order (equivalently, the Bely\u\i\/ functions $X\to X/G_i\cong {\mathbb P}^1({\mathbb C})$ have the same degree), and in particular where these groups are isomorphic. Regular dessins on the sphere and the torus are well known and easily described, so to avoid trivial cases and special arguments we will restrict attention to dessins of genus $g\ge 2$, for which $\mathbb U$ is the hyperbolic plane $\mathbb H$.

The most obvious cause of this phenomenon is a non-normal inclusion $\Delta<\Delta^*$ between triangle groups: if a subgroup $K$ of $\Delta$ is normal in $\Delta^*$ then the conjugates $\Delta_i$ of $\Delta$ in $\Delta^*$ induce isomorphic but distinct regular dessins ${\cal H}_i$ on $X={\mathbb H}/K$, with $G_i={\rm Aut}\,{\cal H}_i\cong\Delta/K$. Girondo~\cite{Gir} has studied this situation in detail, so we will consider it only briefly in Section~4, giving some examples and showing that the number of such dessins ${\cal H}_i$ appearing on $X$ is at most $10$, attained only when $\Delta$ and $\Delta^*$ have types $(3,8,8)$ and $(2,3,8)$ respectively.

A less common example of this phenomenon occurs when triangle groups $\Delta_1$ and $\Delta_2$ of the same type are non-conjugate subgroups of a triangle group $\Delta^*$. Girondo and Wolfart~\cite{GW} have shown that in this case the groups $\Delta_i$ and $\Delta^*$ have types $(n,2n,2n)$ and $(2,2n,2n)$ for some $n$, with $|\Delta^*:\Delta_i|=2$. A surface group $K$ which is normal in both groups $\Delta_i$ corresponds to a pair of regular dessins ${\cal H}_i$ of the same type, with possibly non-isomorphic automorphism groups $G_i\cong\Delta_i/K$ of the same order, on the same surface $X={\mathbb H}/K$. Girondo and Wolfart give such an example, of genus $4$, for $n=3$, with $G_1\not\cong G_2$; in Section~5 we will show that this is the first of an infinite series of examples of genus $(n-1)^2$ for $n\geq 3$, lying on certain quotients of Fermat curves. We will also show how to construct infinite families of examples for which $G_1\cong G_2$, even though ${\cal H}_1\not\cong{\cal H}_2$.

The only other cases in which we can have non-isomorphic dessins ${\cal H}_i$ on the same surface, with automorphism groups $G_i$ of the same order, are when $\Delta_1$ and $\Delta_2$ have types $(2n,2n,2n)$ and $(n,4n,4n)$ for some $n$, and $\Delta^*$ has type $(2,2n,4n)$ or $(2,3,4n)$. Examples of these two cases are constructed in Sections~6 and 7.

Even in the restricted context of these four cases, the residual finiteness of triangle groups and the abundance of finite groups together make it unrealistic to expect complete classifications of the relevant pairs of dessins ${\cal H}_i$. Instead, we will give group-theoretic constructions of infinite families of examples, hoping that some of these methods may also be useful elsewhere. To simplify the exposition we have restricted attention to elementary examples of groups, such a symmetric groups and projective groups over prime fields, but in principle these methods can be applied to much wider classes of groups.

The author is grateful to J\"urgen Wolfart for some useful suggestions about Example~4, and to the organisers of the workshop `Groups and Languages' in honour of Antonio Mach\`\i\/, at the Universit\`a di Roma `La Sapienza', for the opportunity to present some of these ideas in the presence of the founders of the theory of hypermaps.

\section{Hypermaps, dessins and triangle groups}

A {\em hypermap\/} $\cal H$ (always assumed to be finite and oriented) is defined to be an ordered pair of permutations $x$ and $y$ (denoted by $\sigma$ and $\alpha$ by Cori and Mach\`\i~\cite{Cor, CoMa}) generating a finite group $G$ acting transitiviely on a set $\Omega$. One can identify $\Omega$ with the set of edges of a bipartite graph, embedded as a map on a compact oriented surface, so that the cycles of $x$ and $y$, known as {\em hypervertices\/} and {\em hyperedges}, correspond to the rotation of edges around the white and black vertices of the graph; this map is the {\em Walsh map\/} $W({\cal H})$ of $\cal H$, with the cycles of $z:=(xy)^{-1}$, known as {\em hyperfaces}, in bijective correspondence with the faces of $W({\cal H})$~\cite{Wal}. The {\em genus\/} of $\cal H$ is that of its underlying surface.

Since $xyz=1$, the Riemann Existence Theorem shows that there is a covering $\beta:X\to{\mathbb P}^1({\mathbb C})={\mathbb C}\cup\{\infty\}$ of the Riemann sphere by a compact Riemann surface $X$, unramified outside $\{0, 1, \infty\}$, with monodromy permutations $x, y$ and $z$ at these points. We call $\beta$ a {\em Bely\u\i\/ function}. The inverse image of the unit interval $[0,1]$ is a bipartite map on $X$, isomorphic to $W({\cal H})$, with white and black vertices lying over $0$ and $1$, and face-centres over $\infty$. As a compact Riemann surface, $X$ can be regarded as a projective algebraic curve defined over $\mathbb C$. Bely\u\i's Theorem~\cite{Bel} shows that the curves $X$ arising in this way from hypermaps are those defined over the field $\overline{\mathbb Q}$ of algebraic numbers. Following Grothendieck~\cite{Gro}, we will use the term {\em dessin\/} to denote a hypermap, represented in the above way as a bipartite map and carrying this analytic and algebraic structure on its underlying surface. This extra structure provides some of the motivation for what follows, but from a combinatorial point of view much of it can be ignored, and our results can be regarded as referring to hypermaps in their original sense.

The automorphism group ${\rm Aut}\,{\cal H}$ of a dessin $\cal H$, defined as the centraliser of $G$ in the symmetric group on $\Omega$, is the orientation- and colour-preserving automorphism group of the bipartite map $W({\cal H})$. By permuting the names hypervertices, hyperedges and hyperfaces (equivalently the critical values $0, 1$ and $\infty$ in ${\mathbb P}^1({\mathbb C})$), we obtain new dessins called the {\em associates\/} of $\cal H$, with the same automorphisms and underlying surface. We call $\cal H$ a {\em regular dessin\/} (or an {\em orientably regular hypermap\/}) if ${\rm Aut}\,{\cal H}$ is transitive on $\Omega$; equivalently $G$ is a regular permutation group, in which case ${\rm Aut}\,{\cal H}\cong G$, both acting regularly, and $\beta$ is a regular covering, with covering group $G$. For the rest of this paper, $\cal H$ or ${\cal H}_i$ will always denote a regular dessin, represented as its Walsh map.

If $\cal H$ is regular, and $x, y$ and $z$ have orders $l, m$ and $n$, then the black and white vertices all have valencies $l$ and $m$, while the faces are all $2n$-gons; we say that $\cal H$ has {\em type\/} $(l, m, n)$. A dessin of type $(l,2,n)$ can be regarded as an $l$-valent map, by ignoring the black vertices of valency $2$ in $W({\cal H})$, so that all faces are $n$-gons; following Coxeter and Moser~\cite{CoMo}, we say that such a map has {\em type\/} $\{n, l\}$. 

Given integers $l, m, n\ge 2$, the triangle group $\Delta=\Delta(l, m, n)$ has a presentation
\[\Delta(l, m, n)=\langle u, v, w \mid u^l=v^m=w^n=uvw=1\rangle.\]
There is an action of $\Delta$ as a group of automorphisms of a simply-connected Riemann surface $\mathbb U$, namely the Riemann sphere ${\mathbb P}^1({\mathbb C})$, the complex plane $\mathbb C$, or the hyperbolic plane $\mathbb H$, as $l^{-1}+m^{-1}+n^{-1}>1$, $=1$ or $<1$; the {\em canonical generators\/} $u, v$ and $w$ are rotations through $2\pi/l$, $2\pi/m$ and $2\pi/n$ about the vertices of a triangle with internal angles $\pi/l$, $\pi/m$ and $\pi/n$. We will be concerned with the last case, when $\Delta$ is a Fuchsian group.

If $\cal H$ is regular, of type $(l, m, n)$, with automorphism group $G$, there is an epimorphism $\Delta=\Delta(l, m, n)\to G$, giving a transitive permutation representation of $\Delta$ on $\Omega$. The kernel $K$ is a torsion-free normal subgroup of finite index in $\Delta$, with $\Delta/K\cong G$, and conversely every such subgroup $K$ of $\Delta$ arises in this way for some regular dessin $\cal H$. The associated Riemann surface $X$ is isomorphic to ${\mathbb U}/K$, with $K$ isomorphic to the fundamental group $\pi_1X$ of $X$, and the Bely\u\i\/ function $\beta$ corresponding to the projection ${\mathbb U}/K\to{\mathbb U}/\Delta\cong{\mathbb P}^1({\mathbb C})$. We say that $X$ is {\em uniformised\/} by the {\em surface group\/} $K$. Isomorphism of Riemann surfaces $X$ is equivalent to conjugacy of the corresponding surface groups $K$ in ${\rm Aut}\,{\mathbb U}$ ($=PSL_2({\mathbb R})$ when ${\mathbb U}={\mathbb H}$). Taking an associate of $\cal H$ corresponds to permuting the periods of $\Delta$ by changing the canonical generators: for instance, transposing vertex-colours corresponds to regarding $\Delta$ as $\Delta(m, l, n)$ with generators $v, u^v$ and $w$.

\section{Multiply quasiplatonic surfaces and triangle group inclusions}

From now on we will assume that

\smallskip
\centerline{\em ${\cal H}_1$ and ${\cal H}_2$ are regular dessins on the same Riemann surface $X$ of genus $g\ge 2$.}
\smallskip
\noindent To avoid trivial cases, we will also assume that ${\cal H}_1$ and ${\cal H}_2$ are not associates of each other. This is equivalent to $X$ being uniformised by a surface group $K\le PSL_2({\mathbb R})$, which is a normal subgroup of distinct triangle groups $\Delta_i$ of the same types as the dessins ${\cal H}_i$, with each $\Delta_i/K\cong G_i:={\rm Aut}\,{\cal H}_i$. Since each $\Delta_i$ normalises $K$ it is contained in the normaliser $N(K)$ of $K$ in $PSL_2({\mathbb R})$; since this is a Fuchsian group containing a triangle group, it follows from results of Singerman~\cite{Sin2} that $N(K)$ must also be a triangle group $\Delta^*$, and that the possibilities for these triangle group inclusions are all known.

Excluding the rather trivial cases of cyclic and dihedral groups, which are not relevant here, Singerman's paper~\cite{Sin2} lists the normal and non-normal inclusions between triangle groups. We will use the notations $\Delta\triangleleft_i\Delta^*$ and $\Delta<_i\Delta^*$ to denote that $\Delta$ is a normal or non-normal subgroup of index $i$ in $\Delta^*$.

The normal inclusions between hyperbolic triangle groups have the forms

\smallskip
\noindent (a) $\Delta(s,s,t)\triangleleft_2\Delta(2,s,2t)$ where $(s-2)(t-1)>2$,
with quotient group $C_2$,

\smallskip
\noindent (b) $\Delta(t,t,t)\triangleleft_3\Delta(3,3,t)$ where $t>3$, with quotient group $C_3$,

\smallskip
\noindent (c) $\Delta(t,t,t)\triangleleft_6\Delta(2,3,2t)$ where $t>3$, with quotient group $S_3$.
\smallskip

The non-normal inclusions between hyperbolic triangle groups have the forms

\smallskip
\noindent (A) $\Delta(7,7,7)<_{24}\Delta(2,3,7)$,\quad\;
(B) $\Delta(2,7,7)<_9\Delta(2,3,7)$,\quad\,
 (C) $\Delta(3,3,7)<_8\Delta(2,3,7)$,
 
\smallskip
\noindent (D) $\Delta(4,8,8)<_{12}\Delta(2,3,8)$, \quad\,
 (E) $\Delta(3,8,8)<_{10}\Delta(2,3,8)$, \quad
 (F) $\Delta(9,9,9)<_{12}\Delta(2,3,9)$,
 
\smallskip
\noindent (G) $\Delta(4,4,5)<_6\Delta(2,4,5)$, \qquad \qquad \qquad \qquad
 (H) $\Delta(n,4n,4n)<_6\Delta(2,3,4n)$ where $n\ge 2$,
 
\smallskip
\noindent (I) $\Delta(n,2n,2n)<_4\Delta(2,4,2n)$ where $n\ge 3$,\quad
(J) $\Delta(3,n,3n)<_4\Delta(2,3,3n)$ where $n\ge 3$, 
 
\smallskip
\noindent (K) $\Delta(2,n,2n)<_3\Delta(2,3,2n)$ where $n\ge 4$.

\medskip

Recent papers~\cite{Gir, GTW, GW, Sin3, SS} by Girondo, Singerman, Syddall, Torres and Wolfart give useful information about these inclusions and their consequences for dessins. From now on we will make the extra assumption that

\smallskip
\centerline{\em the automorphism groups $G_i$ of the two regular dessins ${\cal H}_i$ have the same order.}
\smallskip

\noindent The corresponding triangle groups $\Delta_1$ and $\Delta_2$ then have the same index $|\Delta^*:K|/|G_i|$ in $\Delta^*:=N(K)$. Inspection of the above lists shows that there are just four possibilities:

\begin{prop}
Let $\Delta_1$ and $\Delta_2$ be two hyperbolic triangle groups of the same index in a triangle group $\Delta^*$. Then one of the following holds:
\begin{enumerate}
\item $\Delta_1$ and $\Delta_2$ have the same type, and are conjugate in $\Delta^*$;
\item $\Delta_1$ and $\Delta_2$ have the same type and are not conjugate in $\Delta^*$;
\item $\Delta_1$ and $\Delta_2$ have types $(2n, 2n, 2n)$ and $(n, 4n, 4n)$ for some $n$, with $\Delta^*=\Delta(2, 2n, 4n)$;
\item $\Delta_1$ and $\Delta_2$ have types $(2n, 2n, 2n)$ and $(n, 4n, 4n)$ for some $n$, with $\Delta^*=\Delta(2, 3, 4n)$.
\end{enumerate}
\end{prop}

Girondo and Wolfart~\cite{GW} have shown that in Case~2, $\Delta_1$ and $\Delta_2$ have the same type $(n, 2n, 2n)$, and are subgroups of index $2$ in $\Delta^*=\Delta(2, 2n, 2n)$, so both inclusions are of type~(a) with $s=2n$ and $t=n$. In Case~3 the inclusions are again of type (a), one of the form $\Delta(2n, 2n, 2n)\triangleleft_2\Delta(2, 2n, 4n)$ with $s=t=2n$, and the other of the form $\Delta(4n, 4n, n)\triangleleft_2\Delta(2,4n,2n)$ with $s=4n$ and $t=n$. In Case~4 the inclusions are of types (c) and (H) with $t=2n$. We will consider these four cases in the following sections.

\section{Case 1}

It has been shown by Girondo and Wolfart in~\cite[Theorem 13]{GW} that if $\Delta_1$ and $\Delta_2$ have the same type, then except in one special case, which we shall consider as Case~2, they are conjugate in $\Delta^*$. In this situation $X$ has an automorphism transforming ${\cal H}_1$ to a distinct but isomorphic dessin ${\cal H}_2$. This phenomenon has been considered in depth by Girondo in~\cite{Gir}, so here we will consider just one typical inclusion, namely~(C), where the subgroups $\Delta_i$ are conjugates of $\Delta=\Delta(3, 3, 7)$ in $\Delta^*=\Delta(2, 3, 7)$.

\medskip

\noindent{\bf Example 1.} The core (intersection of conjugates) of $\Delta$ in $\Delta^*$ is the surface group $K$ of genus $3$ uniformising Klein's quartic curve $X$, given as a projective algebraic curve by the equation
\[x^3y+y^3z+z^3x=0.\]
This Riemann surface has automorphism group
\[G={\rm Aut}\,X\cong \Delta^*/K\cong L_2(7)\cong L_3(2),\]
where $L_n(q):=PSL_n({\mathbb F}_q)$. If we regard $\Delta^*$ as $\Delta(3, 2, 7)$, the regular dessin ${\cal H}^*$ corresponding to its normal subgroup $K$ is a trivalent map tessellating $X$ by $21$ $7$-gons, with ${\rm Aut}\,{\cal H^*}\cong G$; this is the dual of the map R3.1 in Conder's computer-generated list of regular maps and hypermaps~\cite{Con2}. There is a single conjugacy class of eight subgroups $\Delta=\Delta_i\;(i=1,\ldots, 8)$ of index $8$ in $\Delta^*$, corresponding to the point-stabilisers in the natural action of $G$ as $L_2(7)$ on ${\mathbb P}^1(7)$ (equivalently its action by conjugation on its Sylow $7$-subgroups). Thus there are eight regular dessins ${\cal H}_i$ of type $(3, 3, 7)$ on $X$, each with
\[{\rm Aut}\,{\cal H}_i\cong\Delta/K\cong C_7:C_3,\]
a nonabelian group of order $21$. These dessins are permuted transitively by $G$, the action being equivalent to that of $G$ as $L_2(7)$ on ${\mathbb P}^1(7)$, so they are mutually isomorphic, appearing as CH3.1 in~\cite{Con2}. Each ${\cal H}_i$, represented as a bipartite map, has seven white and seven black vertices, each of valency $3$, with both sets permuted transitively by a subgroup of order $7$; the points of $X$ appearing as these vertices, for $i=1,\ldots, 8$, are the vertices of the map ${\cal H}^*$, each appearing in two dessins ${\cal H}_i$. Each dessin ${\cal H}_i$ has three faces, each a $14$-gon, incident with every vertex of ${\cal H}_i$; the centres of these faces are also the face-centres of ${\cal H}^*$, each appearing in one dessin ${\cal H}_i$. The $24$ face-centres of ${\cal H}^*$ are thus partitioned into eight sets of three, each set appearing as the face-centres in ${\cal H}_i$ for some $i$, with the three elements of the set fixed by a subgroup $C_7$ of $G$ fixing a particular point in ${\mathbb P}^1(7)$.

\medskip

\noindent{\bf Example 2.} In any other example based on inclusion~(C) the underlying Riemann surface is a covering of the surface $X$ in Example~1, corresponding to a surface group contained in $K$. For instance, following Macbeath's construction of an infinite sequence of Hurwitz groups~\cite{MacHur}, for any integer $m\geq 2$ let $L=K'K^m$ be the characteristic subgroup of $K$ generated by its commutators and $m$th powers. This is a normal subgroup of $\Delta^*$, contained in each of the eight subgroups $\Delta_i$. Since $K$ is a surface group of genus $g=3$, the quotient $K/L$ is isomorphic to $C_m^{2g}=C_m^6$. As a subgroup of finite index in a surface group, $L$ is also a surface group. The corresponding surface $Y={\mathbb H}/L$ has genus $|K:L|(g-1)+1=2m^6+1$ since the covering $Y\to X$ is unbranched. The normal inclusions $L\triangleleft\Delta_i$ correspond to eight isomorphic regular dessins $\tilde{\cal H}_i$ on $Y$, each an $m^6$-sheeted regular covering of ${\cal H}_i$. These are permuted transitively by $\tilde G={\rm Aut}\,Y\cong \Delta^*/L$, an extension of $K/L$ by $G$. The kernel of this action is $K/L$, and the induced action is again equivalent to that of $G$ on ${\mathbb P}^1(7)$. 

\medskip

\noindent{\bf Example 3.} In Examples~1 and 2, $K$ is identified with $\pi_1X$, its abelianisation $K^{\rm ab}=K/K'$ is identified with the first integer homology group $H_1(X;{\mathbb Z})$, and $K/L$ is identified with the mod~$(m)$ homology group $H_1(X;{\mathbb Z}_m)=H_1(X;{\mathbb Z})\otimes_{\mathbb Z}{\mathbb Z}_m$. The natural action of $G$ on these homology groups, induced by its action on $X$, is equivalent to the action of $\Delta^*/K$ by conjugation on $K/K'$ and on $K/L$. Any $G$-invariant submodule $M/L$ of $K/L$, corresponding to a normal subgroup $M$ of $\Delta^*$ lying between $K$ and $L$, will give rise to further coverings of the dessins ${\cal H}_i$ on the surface $Z={\mathbb H}/M$. For instance, if we take $m=2$ then $K/L=H_1(X;{\mathbb Z}_2)$ is a direct sum of two $3$-dimensional irreducible $G$-submodules, corresponding to the Brauer characters $\varphi_2$ and $\varphi_3$ in~\cite{JLPW}. These characters are complex conjugates of each other, so $X$ has a chiral pair of $8$-sheeted regular covering surfaces $Z$, each carrying eight isomorphic regular dessins of type $(3, 3, 7)$ and genus $17$ (CH17.1 in~\cite{Con2}) permuted transitively by ${\rm Aut}\,Z\cong\Delta^*/M$.

\medskip

There are similar examples based on inclusion~(C) for other values of $m$ (see, for instance, Cohen's construction of Hurwitz groups as abelian coverings of $L_2(7)$ in~\cite{Coh}), and also on the other inclusions in Singerman's list.

For each of these inclusions, the number of mutually isomorphic dessins ${\cal H}_i$ on $X$ is equal to the number of conjugates of $\Delta$ in $\Delta^*$, and this is the index $|\Delta^*:N_{\Delta^*}(\Delta)|$ in $\Delta^*$ of the normaliser $N_{\Delta^*}(\Delta)$ of $\Delta$ in $\Delta^*$. This normaliser can be determined by considering the action of $\Delta^*$ on the cosets of $\Delta$, using ideas described by Singerman in~\cite{Sin1}. In the case of the normal inclusions (a), (b) and (c), of course, the index is $1$. In cases (B), (C), (E), (G), (J) and (K) we have $N_{\Delta^*}(\Delta)=\Delta$, of indices $9$, $8$, $10$, $6$, $4$ and $3$. In case (A) we have $N_{\Delta^*}(\Delta)=\Delta(3,3,7)$, of index $8$ in $\Delta^*$, an inclusion of type (C). In case (D), $N_{\Delta^*}(\Delta)=\Delta(2,8,8)$ of index $6$ and type (H). In case (F), $N_{\Delta^*}(\Delta)=\Delta(3,3,9)$ of index $4$ and type (J). In case (H), $N_{\Delta^*}(\Delta)=\Delta(4n,2,2n)$ of index $3$ and type (K). In case (I), $N_{\Delta^*}(\Delta)=\Delta(2n,2,2n)$ of index $2$ and type (a).

The number of isomorphic dessins ${\cal H}_i$ on $X$ is thus at most $10$, attained in case~(E), where $\Delta^*=\Delta(2,3,8)$ acts on the cosets of $\Delta=\Delta(3,8,8)$ as $PGL_2(9)$ on ${\mathbb P}^1(9)$. The first example of this, with $K$ the core of $\Delta$ in $\Delta^*$, is a surface of genus $16$,; this carries a map of type $\{3, 8\}$, denoted by R16.1 in~\cite{Con2}, with $\Delta^*$ regarded as $\Delta(8,2,3)$.

\section{Case 2}

Girondo and Wolfart~\cite[Theorem 13]{GW} have shown that the only case in which a hyperbolic triangle group $\Delta^*$ contains two non-conjugate triangle groups of the same type is when $\Delta_1$ and $\Delta_2$, both of type $(n, 2n, 2n)$ for some $n\geq 3$, are distinct subgroups of index $2$ in $\Delta^*=\Delta(2, 2n, 2n)$. Specifically, $\Delta_1$ and $\Delta_2$ are the normal closures in $\Delta^*$ of its second and third canonical generators, both of order $2n$; these inclusions are of Singerman's type~(a) with $s=2n$ and $t=n$. Although $\Delta_1$ and $\Delta_2$ are not conjugate in $\Delta^*$, they are conjugate in a triangle group $\Delta^{\dagger}=\Delta(2, 4, 2n)$ which contains $\Delta^*$ as a subgroup of index $2$.

If a surface group $K$ is normal in both $\Delta_1$ and $\Delta_2$, and hence in $\Delta^*$, then by regarding each $\Delta_i$ as $\Delta(2n, 2n, n)$ we obtain dessins ${\cal H}_i$ of type $(2n, 2n, n)$ on the same surface $X={\mathbb H}/K$ of genus
\[g=1+\frac{n-2}{2n}|\Delta_i:K|,\]
with automorphism groups $G_i\cong\Delta_i/K$. These dessins can be represented as bipartite maps, and if we ignore the vertex colours these are simply maps ${\cal M}_i$ of type $\{2n, 2n\}$ on $X$. Since $K$ is normal in $\Delta^*$, which we may regard as $\Delta(2n, 2, 2n)$, these maps are regular dessins, with ${\rm Aut}\,{\cal M}_i\cong\Delta^*/K$. These maps are mutually dual, with conjugation in $\Delta^{\dagger}$ transposing the two subgroups $\Delta_i$ and also the two canonical generators of $\Delta^*$ of order $2n$.

The inclusion of $K$ in $\Delta^{\dagger}$ corresponds to the median map ${\cal M}^{\dagger}={\rm Med}\,{\cal M}_i$ of both maps ${\cal M}_i$, a map of type $\{2n, 4\}$ on $X$. Given any map $\cal M$, the {\em median map\/} ${\rm Med}\,{\cal M}$ is a map on the same surface as $\cal M$; its vertices are the midpoints $m_e$ of the edges $e$ of ${\cal M}$, and if $e$ and $e'$ are consecutive edges of a face $f$ of ${\cal M}_i$ then $m_e$ and $m_{e'}$ are joined by an edge in $f$; thus ${\rm Med}\,{\cal M}$ has valency $4$, each $k$-valent vertex of $\cal M$ lies in a $k$-gonal face of ${\rm Med}\,{\cal M}$, and each $l$-gonal face of $\cal M$ encloses a smaller $l$-gonal face of ${\rm Med}\,{\cal M}$. A map and its dual have the same median map. In our case, $k=l=2n$, so ${\cal M}^{\dagger}$ has type $\{2n, 4\}$. As a dessin, ${\cal M}^{\dagger}$ has automorphism group ${\rm Aut}\,{\cal M}^{\dagger}\cong N_{\Delta^{\dagger}}(K)/K$, isomorphic to $\Delta^{\dagger}/K$ or to $\Delta^*/K$ as $K$ is or is not normal in $\Delta^{\dagger}$, that is, as ${\cal M}_1\cong{\cal M}_2$ or not.

The following example is considered in some detail, for two reasons. Firstly it shows that an important example described by Girondo and Wolfart in~\cite{GW}, where $n=3$, is in fact the first of an infinite family of examples, arising for all $n\geq 3$. Secondly, it is one of the few instances where it is possible to give a completely explicit link between the algebraic geometry of the curve and the combinatorial and group-theoretic aspects of the dessins. As usual in such situations, it is simplest to start with the curve.

\medskip

\noindent{\bf Example 4.} Let $X$ be the Riemann surface corresponding to the affine curve $x^{2n}+y^n=1$, where $n\geq 3$. This has `obvious' automorphisms
\[a:(x,y)\mapsto (\zeta_{2n}x,y)\quad {\rm and}\quad b:(x,y)\mapsto (x,\zeta_ny)
\qquad(\zeta_m:=e^{2\pi i/m}),\]
generating a subgroup $G_1\cong C_{2n}\times C_n$ of ${\rm Aut}\,X$. The regular covering $X\to X/G_1\cong{\mathbb P}^1({\mathbb C})$ is given by the Bely\u\i\/ function
\[\beta:(x,y)\mapsto x^{2n},\]
with critical values $0$ and $\infty$, over each of which there are $n$ points of multiplicity $2n$, and $1$, over which there are $2n$ points of multiplicity $n$. It follows from the Riemann-Hurwitz formula that $X$ has genus $(n-1)^2$.

Since $a$, $b$ and $ab$ have orders $2n$, $n$ and $2n$, there is a smooth epimorphism $\theta:\Delta_1=\Delta(2n, n, 2n)\to G_1$, sending the canonical generators of $\Delta_1$ to $a, b$ and $(ab)^{-1}$. We can identify $X$ with ${\mathbb H}/K$, where $K={\rm ker}\,\theta$ is a normal surface subgroup of $\Delta_1$ with $\Delta_1/K\cong G_1$. Since $G_1$ is abelian we have $K\geq\Delta_1'$, and since $|\Delta_1:K|=2n^2=|\Delta_1:\Delta_1'|$ we have $K=\Delta_1'$. Thus $K$ is a characteristic subgroup of $\Delta_1$, so it is a normal subgroup of $N(\Delta_1)$; by Singerman's results~\cite{Sin2}, $N(\Delta_1)$ is a triangle group $\Delta^*\cong\Delta(2, 2n, 2n)$, containing $\Delta_1$ with index $2$, so $G_1$ has index $2$ in a subgroup $G\cong\Delta^*/K$ of ${\rm Aut}\,X$.

In order to determine this larger group, note that the canonical generator of order $2$ of $\Delta^*$ acts by conjugation on $\Delta$ by transposing its canonical generators of order $2n$. These induce the automorphisms $a$ and $(ab)^{-1}$ of $X$, each of which has $n$ fixed points: those of $(ab)^{-1}$ are the points $P_j\;(j\in{\mathbb Z}_n)$ above $\infty$, where $x, y\to\infty$ with $x^2/y$ approaching a specific $n$th root $\zeta_{2n}^j$ of $-1$, for odd $j\in{\mathbb Z}_{2n}$, while those of $a$ are the points $Q_j=(0,\zeta_n^j)$ above $0$, for $j\in{\mathbb Z}_n$. This suggests adjoining to $G_1$ an automorphism of the form $c:(x,y)\mapsto(\lambda/x, r(x,y))$ for some $\lambda\in{\mathbb C}$ and some rational function $r(x,y)$ satisfying
\[\frac{\lambda^{2n}}{x^{2n}}+r^n=1\]
whenever $x^{2n}+y^n=1$. An obvious choice is to put $\lambda^{2n}=1$ and $r(x,y)=\mu y/x^2$ where $\mu^{n}=-1$. This gives an automorphism $c$ satisfying
\[c^2:(x,y)\mapsto\Bigl(x,\frac{\mu^2y}{\lambda^2}\Bigr),\]
so if we take $\lambda=\mu=\zeta_{2n}$ we obtain an involution
\[c:(x,y)\mapsto\Bigl(\frac{\zeta_{2n}}{x},\frac{\zeta_{2n}y}{x^2}\Bigr).\]
This automorphism, which is clearly not in $G_1$, commutes with $b$ and satisfies $a^c=a^{-1}b^{-1}$, so the subgroup $G:=\langle a, b, c \rangle$ of ${\rm Aut}\,X$ is a semidirect product of $G_1$ by $\langle c \rangle\cong C_2$. The centre of $G$ is $Z=\langle a^n, b \rangle\cong C_2\times C_n$. There is a subgroup $D=\langle a^2b, c \rangle\cong D_n$ in $G$, and if $n$ is odd then $G=Z\times D\cong C_{2n}\times D_n$ (as noted by Girondo and Wolfart~\cite{GW} in the case $n=3$); however, if $n$ is even then $Z\cap D=\langle a^nb^{n/2}\rangle\cong C_2$ and $|G:ZD|=2$.

We can now determine ${\rm Aut}\,X\cong N(K)/K$. We have $N(K)\geq\Delta^*$, and by Singerman's results~\cite{Sin2} the only Fuchsian group properly containing $\Delta^*$ is a triangle group $\Delta^{\dagger}=N(\Delta^*)\cong\Delta(2,4,2n)$. Acting by conjugation, this transposes the two subgroups $\Delta_1, \Delta_2\cong\Delta(2n, n, 2n)$ of index $2$ in $\Delta^*$. It follows that if $N(K)=\Delta^{\dagger}$ then $\Delta_1/K\cong \Delta_2/K$; one of these is the abelian group $\Delta/K\cong G_1$, so they are both abelian and hence their intersection $(\Delta_1\cap\Delta_2)/K$ is a central subgroup of index $4$ in $\Delta^*/K\cong G$. However, the centre $Z$ of $G$ has index $2n$, so $n\leq 2$, against our choice. Thus $N(K)=\Delta^*$, so ${\rm Aut}\,X=G$.

In order to understand the dessins associated with the inclusions of $K$ in these various triangle groups, it is convenient to replace the Bely\u\i\/ function $\beta(x,y)=x^{2n}$ used above with a Bely\u\i\/ function
\[\beta_1=\frac{1}{1-\beta}: (x,y)\mapsto y^{-n}.\]
This change simply permutes the critical values in a $3$-cycle, so that the points above $0, 1$ and $\infty$ have multiplicities $2n, 2n$ and $n$. We now regard $\Delta_1$ and $\Delta_2$ as having type $(2n, 2n, n)$, so each of the normal inclusions $K\triangleleft\Delta_i$ corresponds to a regular dessin ${\cal H}_i$ of type $(2n, 2n, n)$ on $X$ with $G_i:={\rm Aut}\,{\cal H}_i\cong \Delta_i/K$. As we have seen, $G_1$ is abelian whereas $G_2$ is not, so $G_1\not\cong G_2$ and hence ${\cal H}_1\not\cong{\cal H}_2$.

If we represent ${\cal H}_1$ as a bipartite map on $X$, then there are $n$ white and $n$ black vertices, all of valency $2n$; these are respectively at the points $P_j$ and $Q_j$ with $j\in{\mathbb Z}_n$, where $\beta_1=0$ or $1$. The edges consist of the points where $0\leq \beta_1\leq 1$, that is, $y^n\geq 1$ and hence $x^{2n}\leq 0$. Thus there are $2n^2$ edges $e_{j,k}$, for odd $j\in{\mathbb Z}_{4n}$ and any $k\in{\mathbb Z}_n$, along which $x=r\zeta_{4n}^j$ with real $r
\geq 0$, and $y=s\zeta_n^k$ with real $s=(r^{2n}+1)^{1/n}\geq 1$. Such an edge joins $P_{j-2k}$ to $Q_k$, so there are two edges $e_{j,k}$ and $e_{j+2n,k}$ joining each pair of white and black vertices. The embedded graph is therefore $2K_{n,n}$, formed by doubling the edges of the complete bipartite graph $K_{n,n}$. The $2n$ faces of ${\cal H}_1$ are $2n$-gons with centres at the points $R_j=(\zeta_{2n}^j, 0)$, $j\in\mathbb {Z}_{2n}$. The automorphisms $ab$ and $a$ fix all the white and black vertices respectively, while $b$ fixes all the face-centres.

By ignoring the vertex colours, we can regard ${\cal H}_1$ as a map ${\cal M}_1$ of type $\{2n, 2n\}$ on $X$. Since $K$ is normal in $\Delta^*$, now regarded as $\Delta(2n, 2, 2n)$, this map is a regular dessin, with ${\rm Aut}\,{\cal M}_1 = G = {\rm Aut}\,X$. The canonical generators $u, v$ and $w$ of $\Delta^*$, of orders $2n, 2$ and $2n$, are mapped to the generators $a$, $c$ and $(ac)^{-1}=abc$ of $G$ respectively.
The Bely\u\i\/ function corresponding to this inclusion is
\[\beta^*=4\beta_1(1-\beta_1):(x,y)\mapsto -4\Bigl(\frac{x}{y}\Bigr)^{2n}.\]
The involution $c\in G$ fixes the $2n$ points $(x,y)\in X$ where $x^2=\zeta_{2n}$ and $y^n=2$. These are the midpoints of the edges $e_{1,k}$ and $e_{2n+1,k}$, which are reversed by $c$.

The standard embedding ${\cal S}_{2n}$ of the complete bipartite graph $K_{2n,2n}$ is the dessin on the Fermat curve $x^{2n}+y^{2n}=1$ corresponding to the Bely\u\i\/ function $(x,y)\mapsto x^{2n}$ (see~\cite{Jon10}). The map ${\cal M}_1$ is the quotient of this by the subgroup generated by the automorphism $(x,y)\mapsto(x,-y)$; the surface covering, given by $(x,y)\mapsto (x,y^2)$, is branched over the face-centres $R_j$. One can also realise ${\cal M}_1$ as a double covering of the standard embedding ${\cal S}_n$ of $K_{n,n}$ on the Fermat curve $x^n+y^n=1$, branched over all $2n$ of its vertices.

The subgroup $\Delta_2\cong\Delta(2n, 2n, n)$ of $\Delta^*$ has canonical generators $w^v$, $w$ and $u^2$. Under the epimorphism $\Delta^*\to \Delta^*/K=G$, these are mapped to $a^{-1}c$, $abc$ and $a^2$, generating a subgroup $G_2$ of index $2$ in $G$, with $G_2={\rm Aut}\,{\cal H}_2$ for some regular dessin ${\cal H}_2$ of type $(2n, 2n, n)$ on $X$. Both $G_1$ and $G_2$ contain $\langle a^2,b\rangle\cong C_n\times C_n$ as a subgroup of index $2$. Although $G_1$ is abelian, $G_2$ is not, since $(a^2)^{abc}=a^{-2}b^{-2}$.

The generators
\[a^{-1}c:(x,y)\mapsto\Bigl(\frac{\zeta_n}{x}, \frac{\zeta_{2n}^3y}{x^2}\Bigr)
\quad{\rm and}\quad
abc:(x,y)\mapsto\Bigl(\frac{1}{x}, \frac{\zeta_{2n}y}{x^2}\Bigr)\]
of $G_2$ have fixed points $R_1, R_{n+1}$ and $R_0, R_n$ on $X$ respectively. Since they act on the points $R_j$ by $R_j\mapsto R_{2-j}$ and $R_j\mapsto R_{-j}$, this group has two orbits of length $n$ on these points, consisting of those $R_j$ with $j$ even or odd. Thus the vertices of ${\cal H}_2$ are the face-centres of ${\cal H}_1$, with white or black colours corresponding to a $2$-face-colouring of ${\cal H}_1$. The fixed points of the third generator $a^2$ are those of $a$, namely the $n$ points $Q_j$; these, together with the $n$ points $P_j$, form an orbit of $H_2$, so the face-centres of ${\cal H}_2$ are the white and black vertices of ${\cal H}_1$. This hypermap ${\cal H}_2$ corresponds to the Bely\u\i\/ function
\[\beta_2:(x,y)\mapsto -\frac{(1-x^n)^2}{4x^n}.\]
Thus the critical points of $\beta_2$ are the same as those for $\beta_1$, namely the points $P_j$, $Q_j$ and $R_j$. However, they are partitioned into three sets, lying over the three critical values, in a different way, with the points $P_j$ and $Q_j$ all lying over $\infty$, while half the points $R_j$ lie over $0$ and half lie over $1$.

In order to identify the graph embedded by ${\cal H}_2$, we have
\[x^n=1-2\beta_2\pm 2\sqrt{\beta_2(\beta_2-1)},\]
with $0\leq\beta_2\leq 1$ along the edges of ${\cal H}_2$. Writing $x^n=u+iv$ with $u, v\in{\mathbb R}$, we see that as $\beta_2$ increases from $0$ to $1$ along the unit interval, the two branches of $x^n$ travel along the upper and lower halves of the unit circle $u^2+v^2=1$ from $x^n=1$, via $x^n=\pm i$ when $\beta_2=1/2$, to $x^n=-1$. Taking $n$th roots, we see that along each edge of ${\cal H}_2$, $x$ follows the unit circle from a white vertex $R_j=(\zeta_{2n}^j,0)$ with $j$ even to a black vertex $R_{j\pm 1}=(\zeta_{2n}^{j\pm 1},0)$. As this happens, $x^{2n}$ goes once around the unit circle in the positive or negative direction, starting and finishing at $1$, so $y^n=1-x^{2n}$ goes round the circle $|y^n-1|=1$ in the same direction, starting and finishing at $0$. Taking $n$th roots, we obtain $n$ closed circuits for $y$, all starting and finishing at $0$. Thus each vertex $R_j$ is joined by $n$ edges to each of the vertices $R_{j\pm 1}$, so the embedded graph is $nC_{2n}$, formed by replacing each of the edges of a $2n$-cycle $C_{2n}$ with $n$ edges connecting the same pair of vertices. The automorphism $a$ of $X$ permutes the $2n$ vertices cyclically, while $b$ fixes the vertices and permutes each of these sets of $n$ edges cyclically.

As in the case of ${\cal M}_1$, if we ignore the vertex-colours of ${\cal H}_2$ we obtain a map ${\cal M}_2$ of type $\{2n, 2n\}$ which is a regular dessin on $X$, with ${\rm Aut}\,{\cal M}_2\cong G\cong\Delta^*/K$. This is the dual of ${\cal M}_1$, with the generators $u$ and $w$ of $\Delta^*$ interchanging their roles in defining vertices and faces. The edge-centres are the same as those for ${\cal M}_1$, namely the $2n$ fixed points of $c$. This map ${\cal M}_2$ can be obtained as an $n$-sheeted covering of the map $\{2n,2\}$, a spherical embedding of a $2n$-circuit $C_{2n}$, branched over its vertices. The associated Bely\u\i\/ function is
\[4\beta_2(1-\beta_2):(x, y)\mapsto -\frac{1}{4}\Bigl(\frac{y}{x}\Bigr)^{2n}.\]
This is the reciprocal of the Bely\u\i\/ function $\beta^*$ associated with ${\cal M}_1$, corresponding to the fact that the critical values $0$ and $\infty$ have been transposed.

Since they embed non-isomorphic graphs, the maps ${\cal M}_1$ and ${\cal M}_2$ are not isomorphic.  This corresponds to $K$ not being normal in $\Delta^{\dagger}$. The inclusion of $K$ in $\Delta^{\dagger}$ corresponds to the median map ${\cal M}^{\dagger}={\rm Med}\,{\cal M}_i$ of both maps ${\cal M}_i$, a non-regular dessin of type $\{2n, 4\}$ on $X$ with ${\rm Aut}\,{\cal M}^{\dagger}\cong N_{\Delta^{\dagger}}(K)/K=\Delta^*/K\cong G$. The Bely\u\i\/ function associated with ${\cal M}^{\dagger}$ is
\[\beta^{\dagger}:(x,y)\mapsto-\frac{(\beta^*-1)^2}{4\beta^*}.\]
(Note that $\beta^{\dagger}$ is invariant under $\beta^*\mapsto 1/\beta^*$.)

The surface $X$ and the Bely\u\i\/ functions considered above are defined over $\mathbb R$ (in fact, over $\mathbb Q$), so ${\cal H}_i$ and ${\cal M}_i$ all admit orientation-reversing automorphisms, that is, they are regular as unoriented hypermaps and maps; in each case, the full automorphism group is obtained by adjoining to $G_1, G_2$ or $G$ an involution which inverts two of the canonical generators. For small $n$ they appear in Conder's lists~\cite{Con2}: for instance, if $n=11$, so that $X$ has genus $100$, then ${\cal M}_1$ and ${\cal M}_2$ are the dual pair of regular maps R100.43.
%while ${\cal M}^{\dagger}$ is R100.5. But this map isn't regular?

\medskip

In Example~4, the dessins ${\cal H}_i$ have non-isomorphic automorphism groups. By contrast, the following construction provides examples in which the automorphism groups of the dessins are isomorphic, even though the dessins themselves are not. 

\begin{lemma} Let $\Delta^*=\Delta(l, m, n)$ with $l, m$ and $n$ all even, so that $\Delta^*$ has three subgroups $\Delta_0$, $\Delta_1$ and $\Delta_2$ of index $2$. Let $G$ be a finite group with a unique subgroup $H$ of index $2$, and suppose that $G$ is a smooth quotient $\Delta^*/L$ of $\Delta^*$, so that $L\leq\Delta_j$ for some $j=0, 1$ or $2$. Then $G$ is also a smooth quotient $\Delta_i/K$ of $\Delta_i$ for each $i\ne j$ in $\{0, 1, 2\}$, where $K=\Delta_i\cap L$ for both $i$.
\end{lemma}

\noindent{\sl Proof.} If we factor out the subgroup of $\Delta^*$ generated by the commutators and the squares, we obtain a Klein four-group $V_4\cong\Delta(2, 2, 2)$, so $\Delta^*$ has exactly three subgroups $\Delta_i\;(i=0, 1, 2)$ of index $2$; these are distinguished from each other by containing one of the three canonical generators of $\Delta^*$, and the squares of the other two. If $G$, with a unique subgroup $H$ of index $2$, is a smooth quotient of $\Delta^*$, then the kernel $L$ is a normal surface subgroup of $\Delta^*$, contained in exactly one of the subgroups of index $2$. Renumbering if necessary, we may assume that $L\le\Delta_0$. Then $L\not\leq\Delta_i$ for each $i=1, 2$, so $K_i:=\Delta_i\cap L$ is a subgroup of index $2$ in $L$. Since $\Delta_1\cap\Delta_0=\Delta_2\cap\Delta_0$ we have $K_1=K_2$, $=K$, say. As an intersection of two normal subgroups of $\Delta^*$, one of them a surface group, $K$ is a normal surface subgroup of $\Delta^*$. Then $\Delta^*/K=\Delta_i/K\times L/K\cong \Delta^*/L\times\Delta^*/\Delta_i\cong G\times C_2$, and in particular $G$ is a smooth quotient $\Delta_i/K$ of $\Delta_i$ for each $i=1, 2$. \hfill$\square$

\medskip

(Here the two distinct but isomorphic subgroups $\Delta_i/K$ ($i\ne j$) of $\Delta^*/K$ correspond to the two copies of $G$ in $G\times C_2$: one is the obvious direct factor, and the other is obtained from this by multiplying the elements of the coset $G\setminus H$ by the involution in $C_2$.)

In particular the group $\Delta^*=\Delta(2, 2n, 2n)$ has three subgroups $\Delta_0$, $\Delta_1$ and $\Delta_2$ of index $2$, distinguished by containing the first, second and third of the canonical generators. Both $\Delta_1$ and $\Delta_2$ are triangle groups of type $(n, 2n, 2n)$, whereas $\Delta_0$ is a quadrilateral group $\Delta(2, 2, n, n)$ with signature $(0; 2, 2, n, n)$. Lemma~5.1 immediately implies the following:

\begin{cor} Let $G$ be a finite group with a unique subgroup $H$ of index $2$. Suppose that $G$ is a smooth quotient $\Delta/L$ of $\Delta^*=\Delta(2, 2n, 2n)$, with $L\leq\Delta_0$. Then $G$ is also a smooth quotient $\Delta_i/K$ of $\Delta_i$ for $i=1, 2$, where $K=\Delta_1\cap L=\Delta_2\cap L$. \hfill$\square$
\end{cor}

In these circumstances, the Riemann surface $X={\mathbb H}/K$ admits regular dessins ${\cal H}_i$ ($i=1, 2$) of type $(n, 2n, 2n)$, with automorphism groups ${\rm Aut}\,{\cal H}_i\cong G$, corresponding to the normal inclusions $K\triangleleft\Delta_i$. By the Riemann-Hurwitz Formula, these dessins have genus
\[g=1+\frac{(n-2)}{2n}|G|.\]

Now $K$ is normal in $\langle\Delta_1,\Delta_2\rangle=\Delta^*$, so $N(K)$ is a Fuchsian group containing $\Delta^*$. By~\cite{Sin2}, the only possibilities are that $N(K)$ is $\Delta^*$ or that it is the maximal triangle group $\Delta^{\dagger}=\Delta(2,4,2n)$. If $N(K)=\Delta^{\dagger}$ then $\Delta_1$ and $\Delta_2$ are conjugate in $N(K)$, so ${\cal H}_1\cong{\cal H}_2$. If, on the other hand, $N(K)=\Delta^*$ then since $\Delta_1$ and $\Delta_2$ are not conjugate in $N(K)$ we have ${\cal H}_1\not\cong{\cal H}_2$. These two possibilities correspond to whether or not $G$, as a quotient of $\Delta^*$, has an automorphism transposing the images $y$ and $z$ of the generators of order $2n$. If $G$ has such an automorphism then $L$ is normal in $\Delta^{\dagger}$ and hence so is $K$, so the two dessins are isomorphic. If $G$ does not have such an automorphism then $L$ and $K$ are not normal in $\Delta^{\dagger}$, and the dessins are associates of each other.

The combinatorial explanation for this is as follows. If we regard $\Delta^*$ as $\Delta(2n,2,2n)$ then the normal inclusion of $L$ in $\Delta^*$ corresponds to a map $\cal M$ of type $\{2n, 2n\}$ which is a regular dessin on a surface $Y={\mathbb H}/L$, with automorphism group ${\rm Aut}\,{\cal M}\cong\Delta^*/L\cong G$. Since $K$ is also normal in $\Delta^*$, and has index $2$ in $L$, it corresponds to a regular dessin $2{\cal M}$: this is a map of type $\{2n, 2n\}$, which is an unbranched double covering of $\cal M$, on a surface $X={\mathbb H}/K$ of genus $g=2g'-1$ where $Y$ has genus $g'$. This map is bipartite and $2$-face-colourable (i.e.~its dual map is also bipartite). It is constructed by taking two copies $v_0$ and $v_1$, coloured white and black, of each vertex $v$ of $\cal M$, with an edge between $v_0$ and $w_1$ whenever $vw$ is an edge of $\cal M$; the cyclic rotation of edges $vw$ around the vertices $v$ of $\cal M$ then determines the cyclic rotations of edges around $v_0$ and $v_1$. This map has orientation-preserving automorphism group $\Delta^*/K\cong G\times C_2$, with the direct factors $G$ and $C_2$ preserving and transposing vertex colours. Since $2{\cal M}$ is bipartite, of type $\{2n, 2n\}$, it is the Walsh map of a regular dessin of type $(2n, 2n, n)$ with automorphism group $G$, and ${\cal H}_1$ is an associate of this, of type $(n, 2n, 2n)$. If we apply the same process to the dual map ${\cal M}'$ of $\cal M$ we obtain a second regular dessin ${\cal H}_2$ on $X$, with the same type and automorphism group. These two dessins are isomorphic if and only if $\cal M$ is self-dual, that is, ${\cal M}\cong{\cal M}'$.

\medskip

\noindent{\bf Example 5.} Taking $n=3$, we define a smooth homomorphism $\theta:\Delta(2, 6, 6)\to S_5$ by sending the canonical generators to
\[x=(12)(34),\quad y=(13)(245)\quad {\rm and}\quad z=(14)(253).\]
The image $G$ is transitive, and hence primitive since the degree is prime; $G$ contains a transposition $y^3$, so $G=S_5$ by~\cite[Theorem~13.3]{Wie}. This group has a unique subgroup $H$ of index $2$, namely $A_5$, which contains $x$, so $L:={\rm ker}\,\theta\leq\Delta_0$. It follows from Corollary~5.2 that $K:=\Delta_1\cap L=\Delta_2\cap L$ is a surface group, normal in both $\Delta_1$ and $\Delta_2$, with each $\Delta_i/K\cong S_5$. The corresponding surface $X$ has genus $21$. The triangle groups $\Delta_i\;(i=1, 2)$ determine regular dessins ${\cal H}_i$ of type $(3, 6, 6)$ on $X$, with ${\rm Aut}\,{\cal H}_i\cong S_5$. Conjugation by the permutation $(25)(34)$ transposes $y$ and $z$, so ${\cal H}_1\cong{\cal H}_2$. These dessins are isomorphic to the hypermap RPH21.4 in~\cite{Con2}, with full automorphism group (including orientation-reversing automorphisms) $S_5\times C_2$. They are constructed, as explained above, from the self-dual regular map $\cal M$ of genus $11$ and type $\{6,6\}$ corresponding to $\theta$, denoted by R11.5 in~\cite{Con2}.

\medskip

\noindent{\bf Example 6.} For an example in which the dessins ${\cal H}_i$ resulting from Corollary~5.2 are not isomorphic, let us again take $n=3$, and define a smooth homomorphism $\theta:\Delta(2, 6, 6)\to S_9$ by sending the canonical generators to
\[x=(17)(28)(46)(59),\quad y=(123456)\quad {\rm and}\quad z=(143827)(569).\]
The image $G$ of $\theta$ is transitive, and the stabiliser of $6$ contains $xy^2=(17359)(284)$ and $z^3=(18)(24)(37)$, so $G$ is $2$-transitive and hence primitive. Since $G$ contains a $3$-cycle $(xy^2)^5$ and an odd permutation $y$, it follows from~\cite[Theorem~13.3]{Wie} that $G=S_9$. Since $y$ and $z$ are odd, while $x$ is even, we have $L:={\rm ker}\,\theta\leq\Delta_0$, so Corollary~5.2 shows that $K:=\Delta_1\cap L=\Delta_2\cap L$ is a surface group (of genus $60481$), normal in both $\Delta_1$ and $\Delta_2$, with each $\Delta_i/K\cong S_9$. The triangle groups $\Delta_i\;(i=1, 2)$ determine regular dessins ${\cal H}_i$ of type $(3, 6, 6)$ on the surface $X={\mathbb H}/K$, with ${\rm Aut}\,{\cal H}_i\cong S_9$.

The canonical generators of order $6$ of $\Delta^*=\Delta(2,6, 6)$ are transposed by conjugation in $\Delta(2,4, 6)$, whereas $y$ and $z$, having different cycle-structures, are not transposed by any automorphism of $S_9$. It follows that $L$ is not normal in $\Delta(2, 4, 6)$, and hence neither is $K$, since $L/K$ is the centre (and hence a characteristic subgroup) of $\Delta^*/K\cong S_9\times C_2$. Thus $N(K)=\Delta^*$, so $\Delta_1$ and $\Delta_1$ are not conjugate in $N(K)$, and hence ${\cal H}_1\not\cong{\cal H}_2$.

\medskip

As in Case~1, infinite families of further examples can be constructed as coverings of these, by considering characteristic subgroups of finite index in $K$, or more generally $G$-invariant subgroups of $K$. Similarly, although the above examples have $n=3$, one can also find examples for other values of $n$, as follows.

\medskip

\noindent{\bf Example 7.} Given $n\ge 3$, let $G=PGL_2(p)$ for some prime $p\equiv 2n+1$ mod~$(4n)$; since $2n+1$ and $4n$ are coprime, Dirichlet's Theorem implies that there are infinitely many such primes. This group $G$ has a unique subgroup $H$ of index $2$, namely $L_2(p)$. Let $x$, $y$ and $z$ be the images in $G$ of the matrices
\[
X=\Big(\,\begin{matrix}a&b\cr c&-a\end{matrix}\,\Big),
\quad
Y=\Big(\,\begin{matrix}d&0\cr 0&1\end{matrix}\,\Big)
\quad{\rm and}\quad
Z=(XY)^{-1}=\Big(\,\begin{matrix}ad&b\cr cd&-a\end{matrix}\,\Big)^{-1}
=\Big(\,\begin{matrix}-a/d&-b/d\cr -c&a\end{matrix}\,\Big)
\]
in $GL_2(p)$, where $a^2+bc+1=0$ and $d$ has multiplicative order $2n$ in ${\mathbb F}_p^*$. Then $x$ is an involution in $H$, while $y$ is an element of order $2n$ in $G\setminus H$ (note that $d$ is not a square in ${\mathbb F}_p^*$ since $2n$ does not divide $(p-1)/2$). Now $XY$ has trace $a(d-1)$, so if we choose
\begin{equation}
a=\frac{d+1}{d-1}
\end{equation}
then $XY$ and $Y$ have the same trace and determinant; they therefore have the same eigenvalues, so those of $Z$ are the inverses of those of $Y$, namely $1/d$ and $1$, and hence $z$ also has order $2n$. If $n\geq 6$ then it follows from Dickson's classification of the subgroups of $L_2(q)$~\cite[Ch.~XII]{Dic} that the only maximal subgroups of $H$ which could contain the element $y^2$ of order $n$ are the stabilisers of the elements $0$ and $\infty$ of ${\mathbb P}^1(p)$ and the dihedral group of order $p-1$ leaving $\{0, \infty\}$ invariant. However, $x$ is not an element of any of these maximal subgroups, since $abc\neq 0$ (otherwise $a=0$, giving $d=-1$, or $a^2=-1$, giving $d^2=-1$, so $d$ has order dividing $4$ and hence $n\le 2$). Thus $x$ and $y^2$ generate $H$ and so $x$ and $y$ generate $G$. Thus $G$ is a smooth quotient of $\Delta^*=\Delta(2, 2n, 2n)$, so by Corollary~5.2 we obtain a pair of regular dessins ${\cal H}_i$ of type $(n, 2n, 2n)$, on the same surface, with ${\rm Aut}\,{\cal H}_i\cong G$.

These two dessins are isomorphic if and only if $G$ has an automorphism transposing $y$ and $z$. Let $u\;(=-b/(a+1))$ and $v\;(=(a-1)/c)$ be the fixed points of $z$ in ${\mathbb P}^1(p)$ corresponding to the eigenspaces of $Z$ with eigenvalues $1/d$ and $1$. Then the involution
\[i:t\mapsto\frac{v(t-u)}{t-v}\]
in $G$ transposes $u$ and $v$ with the fixed points $0$ and $\infty$ of $y$, where $Y$ has eigenvalues $1$ and $d$. Now an element of $G$ with two given fixed points is uniquely determined by the ratio of the corresponding eigenvalues of a matrix representing it. It follows that $i$, acting by conjugation, transposes $y$ and $z$, so ${\cal H}_1\cong{\cal H}_2$.

However, we can often make our choice of $a$ differ from that in $(1)$ (with $a^2\ne -1$ as before), so that $Z$ has eigenvalues $\lambda, \mu=d^{(\pm j-1)/2}$ for some unit $j\not\equiv \pm 1$ mod~$(2n)$. Then $\det Z=\lambda\mu=d^{-1}=\det Y^{-1}$ and $\lambda/\mu=d^j$, so $z$ still has order $2n$. As before, $x$ and $y$ generate $G$, and we obtain two regular dessins ${\cal H}_i$ of type $(n, 2n, 2n)$ on the same surface, with ${\rm Aut}\,{\cal H}_i\cong G$. Since $\lambda/\mu\ne d^{\pm 1}$, no inner automorphism of $G$ can transpose $y$ and $z$. As the automorphism group of $H$, a nonabelian simple group, $G$ is complete. Thus each automorphism of $G$ is inner, and hence cannot transpose $y$ and $z$, so ${\cal H}_1\not\cong{\cal H}_2$.

For example, if $n=8$ and $p=17$ we obtain pairs of regular dessins ${\cal H}_i$ of type $(8,16,16)$ on the same surface of genus $1837$, with automorphism group $PGL_2(17)$. Let us take $d=3$, a primitive root mod~$(17)$. If we choose $a=2$, as in $(1)$, we can put $b=3$ and $c=4$ to obtain
\[
X=\Big(\,\begin{matrix}2&3\cr 4&-2\end{matrix}\,\Big),
\quad
Y=\Big(\,\begin{matrix}3&0\cr 0&1\end{matrix}\,\Big)
\quad{\rm and}\quad
Z=\Big(\,\begin{matrix}5&-1\cr -4&2\end{matrix}\,\Big).
\]
Then ${\cal H}_1\cong{\cal H}_2$ since $Y$ and $Z$ have eigenvalues $3, 1$ and $1, 6$ with the same ratio $3^{\pm 1}=6^{\mp 1}$. If, however, we choose $a=-1$, we can put $b=2$ and $c=-1$, giving
\[
X=\Big(\,\begin{matrix}-1&2\cr -1&1\end{matrix}\,\Big),
\quad
Y=\Big(\,\begin{matrix}3&0\cr 0&1\end{matrix}\,\Big)
\quad{\rm and}\quad
Z=\Big(\,\begin{matrix}6&5\cr 1&-1\end{matrix}\,\Big).
\]
In this case $Z$ has eigenvalues $d=3$ and $d^{-2}=2$ with ratio $d^3=-5\ne 3^{\pm 1}$, so ${\cal H}_1\not\cong{\cal H}_2$.

\medskip

The condition $n\ge 6$ was imposed in Example~7 to ensure that $x$ and $y$ generate $G$. However, in some cases they do this even when $n<6$.

\medskip
 
\noindent{\bf Example 8.} Let $n=3$ and $p=7$. Formula~$(3)$ in the Appendix and the character table in~\cite{ATLAS} show that the group $G=PGL_2(7)$ contains $|G|=336$ triples $x, y$ and $z$ of orders $2, 6$ and $6$ with $xyz=1$. By inspection, no proper subgroup contains such a triple, so they all generate $G$. They form a single orbit under ${\rm Aut}\,G=G$, so $\Delta^*=\Delta(2, 6, 6)$ has a unique normal surface subgroup $L$ with $\Delta^*/L\cong G$. Since $x$ is contained in the unique subgroup $H=L_2(7)$ of index $2$ in $G$ we have $L\le\Delta_0$. The two regular dessins ${\cal H}_i$ of type $(3, 6, 6)$ resulting from Corollary~5.2 have genus $57$, and are isomorphic, with ${\rm Aut}\,{\cal H}_i\cong PGL_2(7)$; they appear as RPH57.15 in Conder's list~\cite{Con2}, with full automorphism group $PGL_2(7)\times C_2$. These dessins are constructed, as explained after Corollary~5.2, from the self-dual regular map $\cal M$ of type $\{6,6\}$ and genus $29$ denoted by R29.9 in~\cite{Con2}.

\medskip

We can also use Lemma~5.1 to construct further examples where ${\rm Aut}\,{\cal H}_1\not\cong{\rm Aut}\,{\cal H}_2$, in addition to those in Example~4. The following consequence of Lemma~5.1 is simply Corollary~5.2 with subscripts permuted:

\begin{cor} Let $G$ be a finite group with a unique subgroup $H$ of index $2$. Suppose that $G$ is a smooth quotient $\Delta^*/L$ of $\Delta^*=\Delta(2, 2n, 2n)$, with $L\leq\Delta_1$. Then $G$ is also a smooth quotient $\Delta_i/K$ of $\Delta_i$ for $i=0, 2$, where $K=\Delta_0\cap L=\Delta_2\cap L$. \hfill$\square$
\end{cor}

In this situation $\Delta_1/K=(\Delta_0\cap\Delta_2)/K\times L/K\cong H\times C_2$, so if $G\not\cong H\times C_2$ then the regular dessins ${\cal H}_i$ corresponding to the inclusions $K\le\Delta_i$ for $i=1, 2$ have non-isomorphic automorphism groups $H\times C_2$ and $G$.

\medskip

\noindent{\bf Example 9.} Let $n=6k$ for some integer $k\ge 2$. In the symmetric group $S_d$ of degree $d=2n+1=12k+1$, let
\[x=(1,12k+1)(2,6k+3)(3,2k+4)(4,k+5)(5,6),\quad z=(1,2,\ldots, 12k)(12k+1)\]
and $y=(zx)^{-1}$. Then $x$ and $z$ are odd permutations of orders $2$ and $12k=2n$, while $y$ is an even permutation with cycles
\[(1, 12k+1, 12k, \ldots, 6k+3),\quad (2, 6k+2, 6k+1,\ldots, 2k+4),\]
\[(3, 2k+3, 2k+2,\ldots, k+5),\quad (4, k+4, k+3,\ldots, 6)\quad{\rm and}\quad (5)\]
of lengths $6k, 4k, k, k$ and $1$, so $y$ has order $12k=2n$. Mapping the canonical generators of $\Delta^*$ to $x, y$ and $z$ therefore gives a smooth homomorphism $\theta:\Delta^*\to S_d$. The cycle-structures of $x$ and $z$ show that the image $G=\langle x, z\rangle$ of $\theta$ is doubly transitive. A particular case of a theorem of H\"ochsmann~\cite{Hoc} states that if any doubly transitive group $G$ of degree $d$ has an element of $2$-power order moving just $m$ points, then either  $m\ge d/2$, or $G\ge A_d$, or $G$ is the symplectic group $PSp_{2m}(2)$ with $d=2^{m-1}(2^m-1)$ for some $m>2$. In our case the involution $x$ moves $m=10$ points. Since $d$ is odd, the symplectic case cannot arise; since $d>2m$ we have $G\ge A_d$, and hence $G=S_d$ since $x$ is an odd permutation. The hypotheses of Corollary~5.3 are now satisfied, with $y\in H=A_d$ and $L={\rm ker}\,\theta$. We therefore obtain two regular dessins ${\cal H}_i$ of type $(n,2n,2n)$ on the Riemann surface $X={\mathbb H}/K$ where $K=\Delta_0\cap L=\Delta_2\cap L$. They have non-isomorphic automorphism groups $A_d\times C_2$ and $S_d$.

\section{Case 3}

In Cases~3 and 4 of Proposition~3.1 we have $\Delta_1=\Delta(2n, 2n, 2n)$ and $\Delta_2=\Delta(n, 4n, 4n)$, while $\Delta^*$ is $\Delta(2, 2n, 4n)$ or $\Delta(2, 3, 4n)$ respectively. In either case, we obtain pairs of regular dessins of types $(2n, 2n, 2n)$ and $(n, 4n, 4n)$ on the same Riemann surface $X={\mathbb H}/K$ of genus
\begin{equation}
g=1+\frac{(2n-3)}{4n}|\Delta_i:K|.
\end{equation}
Having different types, these dessins cannot be isomorphic. As we shall see, their automorphism groups may or may not be isomorphic.

In dealing with Case~3, we will simply assume that $N(K)\ge\Delta^*=\Delta(2, 2n, 4n)$, rather than that $N(K)=\Delta^*=\Delta(2, 2n, 4n)$, so that any results obtained here can also be applied in Case~4, where $N(K)$ is the larger group $\Delta(2, 3, 4n)$. In this situation the analogue of Corollary~5.2 is the following immediate consequence of Lemma~5.1, where the quadrilateral group $\Delta_0=\Delta(2, 2, n, 2n)$ with signature $(0; 2, 2, n, 2n)$ is the third subgroup of index $2$ in $\Delta^*$, not equal to $\Delta_1$ or $\Delta_2$:

\begin{cor} Let $G$ be a finite group with a unique subgroup $H$ of index $2$. Suppose that $G$ is a smooth quotient $\Delta^*/L$ of $\Delta^*=\Delta(2, 2n, 4n)$, with $L\leq\Delta_0$. Then $G$ is also a smooth quotient $\Delta_i/K$ of $\Delta_i$ for $i=1, 2$, where $K=\Delta_1\cap L=\Delta_2\cap L$. \hfill$\square$
\end{cor}

The combinatorial explanation for this is similar to that given earlier for Corollary~5.2. As before, the normal inclusions of $L$ and $K$ in $\Delta^*$, now regarded as $\Delta(2n,2,4n)$, correspond to a map $\cal M$ and a bipartite double covering $2{\cal M}$, both of type $\{4n, 2n\}$, and both regular dessins with automorphism groups $G$ and $G\times C_2$. Now $2{\cal M}$ is the Walsh map of a regular dessin ${\cal H}_1$ of type $(2n, 2n, 2n)$ on $X={\mathbb H}/K$ with automorphism group $G$. Applying the same process to the dual map ${\cal M}'$, of type $\{2n, 4n\}$, leads to a regular dessin of type $(4n, 4n, n)$ on $X$, and ${\cal H}_2$ is an associate of this, of type $(n, 4n, 4n)$, also with automorphism group $G$.

\medskip

\noindent{\bf Example 10.} We adapt the construction in Example~9. In the symmetric group $S_d$ of degree $d=4n+1$, where $n\ge 2$, let
\[x=(1,4n+1)(2,2n+3)(3,n+4)(4,5),\quad z=(1,2,\ldots,4n)(4n+1)\]
and $y=(zx)^{-1}$. Then $x$ and $z$ are even and odd permutations of orders $2$ and $4n$, while
\[y^{-1}=zx=(1, 2n+3, 2n+4, \ldots, 4n+1)(2, n+4, n+5, \ldots, 2n+2)(3, 5, 6, \ldots, n+3)(4)\]
has disjoint cycles of lengths $2n, n, n$ and $1$, so $y$ is odd and has order $2n$. We therefore have a smooth homomorphism $\theta:\Delta^*\to S_d$ with $L={\rm ker}\,\theta\le\Delta_0$. As in Example~9, the image $G$ of $\theta$ is doubly transitive, so H\"ochsmann's theorem~\cite{Hoc} shows that $G=S_d$ provided $n\ge 4$. If $n=2$ then $z^3x=(1, 5, 8, 6, 9)(2, 4)$ and $G$ contains $(z^3x)^5=(2, 4)$, while if $n=3$ then $z^3x=(1, 5, 8, 11, 9, 12, 7, 10, 13)(2, 4, 3, 6)$ and $G$ contains $(z^3x)^9=(2, 4, 3, 6)$; in either case it follows as before that $G=S_d$. Corollary~6.1 therefore applies, giving two regular dessins ${\cal H}_i$, of types $(2n,2n,2n)$ and $(n,4n,4n)$, with ${\rm Aut}\,{\cal H}_i\cong G=S_d$, on the Riemann surface $X={\mathbb H}/K$ where $K=\Delta_1\cap L=\Delta_2\cap L$. For instance, if $n=2$, so that $G=S_9$, these dessins have genus $45361$.

\medskip

If $H$ is a smooth quotient $\Delta^*/M$ of $\Delta^*$ with no subgroup of index $2$, then $L:=\Delta_0\cap M$ satisfies the hypotheses of Corollary~6.1, with $G=\Delta^*/L=(\Delta_0/L)\times(M/L)\cong H\times C_2$.

\medskip

\noindent{\bf Example 11.} If $n\ge 2$ then by Dirichlet's Theorem there are infinitely many primes $p\equiv 1$ mod~$(8n)$. For any such $p$, let $x, z$ and $y$ be the images in $H=L_2(p)$ of the matrices
\[
X=\Big(\,\begin{matrix}a&b\cr c&-a\end{matrix}\,\Big),
\quad
Z=\Big(\,\begin{matrix}d&0\cr 0&1/d\end{matrix}\,\Big)
\quad{\rm and}\quad
Y=(ZX)^{-1}=XZ^{-1}=\Big(\,\begin{matrix}a/d&bd\cr c/d&-ad\end{matrix}\,\Big)
\]
in $SL_2(p)$, where $a^2+bc+1=0$ and $d$ has multiplicative order $8n$ in ${\mathbb F}_p^*$. Then $x$ and $z$ have orders $2$ and $4n$, and if we put $a=(d^4+1)/(d-d^3)$ then $Y$ has trace $d^2+d^{-2}$, so $y$ has order $2n$. As in Example~7, Dickson's classification of the maximal subgroups of $L_2(q)$~\cite[Ch.~XII]{Dic} implies that $x$ and $z$ generate $H$, so $H$ is a smooth quotient $\Delta^*/M$ of $\Delta^*$. Since $H$ is simple, it follows from the preceding remark that $G=L_2(p)\times C_2$ satisfies the hypotheses of Corollary~6.1. We therefore obtain two regular dessins ${\cal H}_i$, of types $(2n,2n,2n)$ and $(n,4n,4n)$, with ${\rm Aut}\,{\cal H}_i\cong G$, on the Riemann surface $X={\mathbb H}/K$ where $K=\Delta_1\cap \Delta_2\cap M$. In the smallest case, where $n=2$ and $p=17$, they have genus $613$.

\medskip

It is tempting to try to use the group $G=PGL_2(p)$ here as in Example~7, since it also has $H=L_2(p)$ as its unique subgroup of index $2$. However, this fails when $\Delta^*=\Delta(2, 2n, 4n)$. We need to send the three canonical generators of $\Delta^*$ to elements $x, y$ and $z$ of orders $2, 2n$ and $4n$ in $G$, with $x\in H$ and $y, z\in G\setminus H$. However, if $G$ has elements of order $4n$ then each element $y$ of order $2n$ is the square of an element of order $4n$, so $y\in H$.

In Examples~10 and 11 the dessins ${\cal H}_i$ have isomorphic automorphism groups, but as in Case~2 we can also construct examples where ${\rm Aut}\,{\cal H}_1\not\cong{\rm Aut}\,{\cal H}_2$. The analogue of Corollary~5.3 as the following consequence of Lemma~5.1:

\begin{cor} Let $G$ be a finite group with a unique subgroup $H$ of index $2$. Suppose that $G$ is a smooth quotient $\Delta^*/L$ of $\Delta^*=\Delta(2, 2n, 4n)$, with $L\leq\Delta_1$. Then $G$ is also a smooth quotient $\Delta_i/K$ of $\Delta_i$ for $i=0, 2$, where $K=\Delta_0\cap L=\Delta_2\cap L$. \hfill$\square$
\end{cor}

As in Case~2 it follows that if $G\not\cong H\times C_2$ the regular dessins ${\cal H}_i$ corresponding to the inclusions $K\le\Delta_i$ for $i=1, 2$ have non-isomorphic automorphism groups $H\times C_2$ and $G$.

\medskip

\noindent{\bf Example 12.} We imitate Examples~9 and 10. In $S_d$, where $d=4n+1\ge 9$, let
\[x=(1,4n+1)(2,2n+3)(4,5),\quad z=(1,2,\ldots,4n)(4n+1).\]
and $y=(zx)^{-1}$. Then $x$ and $z$ are odd permutations of orders $2$ and $4n$, while
\[y=(1, 4n+1, 4n, \ldots, 2n+3)(2, 2n+2, 2n+1, \ldots, 3)(4),\]
with cycles of lengths $2n, 2n$ and $1$, is an even permutation of order $2n$. We therefore obtain a smooth homomorphism $\Delta^*\to S_d$. Since $x$ moves six points, H\"ochsmann's theorem~\cite{Hoc} shows that $G:=\langle x, z\rangle=S_d$ provided $n\ge 3$. If $n=2$ then $z^2x=(1, 3, 4, 6, 8, 7, 9)(2, 5)$, so $G$ contains the transposition $(z^2x)^7$ and again $G=S_d$. Corollary~6.2 therefore gives two regular dessins ${\cal H}_1$ and ${\cal H}_2$ of types $(2n,2n,2n)$ and $(n,4n,4n)$ on the same Riemann surface $X={\mathbb H}/K$. We have ${\rm Aut}\,{\cal H}_1\cong A_d\times C_2$ while ${\rm Aut}\,{\cal H}_2\cong S_d\not\cong{\rm Aut}\,{\cal H}_1$. 

\medskip

\noindent{\bf Example 13.} Here we adapt Example~7. Given $n\ge 3$, let $G=PGL_2(p)$ for some prime $p\equiv 4n+1$ mod~$(8n)$; as before, Dirichlet's Theorem implies that there are infinitely many such primes. Let $x$, $y$ and $z$ be the images in $G$ of the matrices
\[
X=\Big(\,\begin{matrix}a&b\cr c&-a\end{matrix}\,\Big),
\quad
Y=\Big(\,\begin{matrix}d&0\cr 0&1/d\end{matrix}\,\Big)
\quad{\rm and}\quad
Z=(XY)^{-1}=\Big(\,\begin{matrix}ad&b/d\cr cd&-a/d\end{matrix}\,\Big)^{-1}
\]
in $GL_2(p)$, where $d$ has multiplicative order $4n$ in ${\mathbb F}_p^*$. Then $x$ and $y$ have orders $2$ and $2n$, with $y\in H=L_2(p)$, and $x\in G\setminus H$ if $a^2+bc$ is a non-square in ${\mathbb F}_p$. If we choose $X$ so that $a^{-2}bc+d+d^{-1}=1$ then $Z^{-1}$ has the same value of ${\rm tr}^2/\det$ as the diagonal matrix with eigenvalues $d$ and $1$, so $z$ has order $4n$. Now $y$ fixes $0$ and $\infty$, while the two points fixed by $z^2$ are disjoint from these provided $bc\ne 0$, so $H=\langle y, z^2\rangle$ and hence $G=\langle y, z\rangle$. Corollary~6.2 then gives two regular dessins of types $(2n,2n,2n)$ and $(n,4n,4n)$ on the same Riemann surface, this time with non-isomorphic automorphism groups $L_2(p)\times C_2$ and $PGL_2(p)$. If $n=3$ and $p=13$, for instance, they have genus $547$.

\section{Case 4}

In Case~4 we have $\Delta_1=\Delta(2n, 2n, 2n)$ and $\Delta_2=\Delta(n, 4n, 4n)$ as in Case~3, but now $\Delta^*:=N(K)=\Delta(2, 3, 4n)$. We have $|\Delta^*:\Delta_i|=6$ for each $i$, with $\Delta_1\triangleleft\Delta^*$ a normal inclusion of type (c), and $\Delta_2<\Delta^*$ a non-normal inclusion of type (H). There is a unique normal subgroup $\Delta_1$ of type $(2n, 2n, 2n)$ in $\Delta^*$, namely the kernel of the natural epimorphism onto $\Delta(2,3,2)\cong S_3$. There is a unique conjugacy class of subgroups $\Delta_2$ of type $(n, 4n, 4n)$ in $\Delta^*$, namely the inverse images of the three subgroups $C_4<S_4$ under the natural epimorphism $\theta:\Delta^*\to \Delta(2,3,4)\cong S_4$ (equivalently, the stabilisers of faces in the action of $\Delta^*$, through $S_4$, by rotations of a cube $\cal C$); here $\Delta_1$ is the inverse image of the normal Klein $4$-group $V_4\triangleleft S_4$. It follows that for each of the three choices of $\Delta_2$ in $\Delta^*$ there is a subgroup $\Delta^{\circ}\cong\Delta(2, 2n, 4n)$ of index $3$ in $\Delta^*$, containing both $\Delta_1$ and $\Delta_2$ with index $2$: these three subgroups $\Delta^{\circ}$ are the inverse images of the Sylow $2$-subgroups of $S_4$ (dihedral groups of order $8$), or equivalently the stabilisers of unordered pairs of opposite faces of $\cal C$.

It follows from these inclusions $\Delta_i\triangleleft_2\Delta^{\circ}<_3\Delta^*$ that Case~4 is just a subcase of the situation considered in Section~6, where the surface group $K$ was normal in $\Delta_1$ and $\Delta_2$, or equivalently in $\langle\Delta_1, \Delta_2\rangle=\Delta^{\circ}\cong\Delta(2, 2n, 4n)$; now it is also normal in the larger group $\Delta^*\cong\Delta(2, 3, 4n)$. The subgroups $K$ we need are therefore the normal surface subgroups of $\Delta^*$ contained in $K_0:=\ker\theta$, corresponding to smooth finite quotients of $\Delta^*$ which map onto $S_4$, or equivalently, to orientably regular maps which cover $\cal C$, branched only over the faces. Any such subgroup $K$ yields regular dessins ${\cal H}_i$ of types $(2n, 2n, 2n)$ and $(n, 4n, 4n)$ on the surface $X={\mathbb H}/K$, with ${\rm Aut}\,{\cal H}_i\cong\Delta_i/K$. In such cases, ${\rm Aut}\,{\cal H}_i$ has $\Delta_i/K_0\cong V_4$ or $C_4$ as an epimorphic image for $i=1$ or $2$; this immediately excludes all the examples constructed in Section~6, where ${\rm Aut}\,{\cal H}_i$ has no normal subgroup of index $4$, so that $K\not\le K_0$.

As in Section~6, the inclusions of $K$ in $\Delta^{\circ}$ correspond to regarding the Walsh maps of ${\cal H}_i\;(i=1, 2)$ as a dual pair of uncoloured maps ${\cal M}_i$ of types $\{4n, 2n\}$ and $\{2n, 4n\}$. The further inclusion of $K$ in $\Delta^*$, now regarded as $\Delta(3, 2, 4n)$, corresponds to the truncation $\cal T$ of ${\cal M}_2$, a map of type $\{4n, 3\}$ on $X$ formed by replacing every vertex (of valency $2n$) of ${\cal M}_2$ with a small $2n$-gon; equivalently, the dual of $\cal T$ is formed by stellating ${\cal M}_1$, joining the centre of each face of ${\cal M}_1$ to all its incident vertices. The normality of $K$ in $\Delta^*$ corresponds to the fact that $\cal T$ and its dual are orientably regular maps, i.e.~regular dessins.

\medskip

\noindent{\bf Example 14.} The simplest examples are those in which $K_0/K$ is cyclic. The cyclic regular coverings $\cal T$ of the cube $\cal C$ were classified by Surowski and the author in~\cite{JSu} (together with those of the other platonic hypermaps). They are Sherk's trivalent maps $\{d\cdot 4, 3\}$ of type $\{4d, 3\}$ and genus $2(d-1)$, described in~\cite{She}, where $n=d:=|K_0:K|$ divides $6$; their orientation-preserving automorphism groups $\Delta^*/K$ are $d$-fold central extensions of $S_4$. If $d=1$ then ${\cal T}={\cal C}$, giving regular dessins ${\cal H}_i$ of types $(2, 2, 2)$ and $(1, 4, 4)$ on the sphere with automorphism groups $V_4$ and $C_4$. If $d=2$ then $\cal T$ is the M\"obius-Kantor map of genus $2$, denoted by $\{4+4, 3\}$ in~\cite[\S8.8, \S8.9, Fig.~3.6c]{CoMo} and by R2.1 in~\cite{Con2}. Here $\Delta^*/K\cong GL_2(3)$, a double covering of $\Delta^*/K_0\cong S_4\cong PGL_2(3)$; the dessin ${\cal H}_1$, denoted by RPH2.4 in~\cite{Con2}, has type $(4, 4, 4)$, and its automorphism group is the quaternion group $Q_8$; the dessin ${\cal H}_2$ has automorphism group $C_8$, and an associate of type $(8, 2, 8)$ appears as the regular map R2.6 in~\cite{Con2}. If $d=3$ then $\cal T$ has genus $4$ and type $\{12, 3\}$, with $\Delta^*/K\cong S_4\times C_3$; the dessins ${\cal H}_i$ have types $(6, 6, 6)$ and $(3, 12, 12)$, and automorphism groups $V_4\times C_3$ and $C_4\times C_3\cong C_{12}$; these maps and hypermaps are R4.1, RPH4.14 and RPH4.8 in~\cite{Con2}. If $d=6$ then $\cal T$ has genus $10$ and type $\{24, 3\}$, with $\Delta^*/K\cong GL_2(3)\times C_3$; the dessins ${\cal H}_i$ have types $(12, 12, 12)$ and $(6, 24, 24)$, and automorphism groups $Q_8\times C_3$ and $C_8\times C_3\cong C_{24}$; these are R10.5, RPH10.38 and RPH10.29 in~\cite{Con2}. In all four cases the underlying surface $X$ is the curve $y^d=x^5-x$, with ${\rm Aut}\,X\cong \Delta^*/K$; the $d$-sheeted covering $(x,y)\mapsto x$ of the sphere ${\mathbb P}^1({\mathbb C})$ is branched over the face-centres $\pm 1, \pm i, 0$ and $\infty$ of $\cal C$.

\medskip

\noindent{\bf Example 15.} Further examples, with $K_0/K$ abelian, can be constructed from quotients
$K_0/K$ of $K_0^{^{\rm ab}}=K_0/K_0'\cong {\mathbb Z}_n^5$, regarded as the homology module \[H_1(Y_0;{\mathbb Z}_n)
=H_1(Y_0;{\mathbb Z})\otimes_{\mathbb Z}{\mathbb Z}_n
=(\pi_1Y_0)^{^{\rm ab}}\otimes_{\mathbb Z}{\mathbb Z}_n\]
over ${\mathbb Z}_n$ for the group $\Delta^*/K_0\cong S_4$, where $Y_0$ is the punctured sphere ${\mathbb P}^1({\mathbb C})\setminus\{0, \pm 1, \pm i, \infty\}$. The induced action of $\Delta^*/K_0$ by conjugation on $K_0^{^{\rm ab}}$ is equivalent to the action of $S_4$ on $H_1(Y_0;{\mathbb Z})$ as the rotation group of $\cal C$. If $n$ is a prime $p>3$, for example, then this $5$-dimensional module is a direct sum of irreducible $S_4$-modules of dimensions $2$ and $3$, so we obtain normal subgroups $K$ of $\Delta^*$ with $K/K_0$ elementary abelian of order $p^e$ for $e=2, 3$ and $5$. If $p^e=25$, for instance, then $\cal T$ is the regular map R36.3 of genus $36$ in~\cite{Con2}.

\medskip

Under the natural epimorphism from the modular group $\Gamma=PSL_2({\mathbb Z})\cong C_2*C_3\cong\Delta(2,3,\infty)$ to $\Delta^*\cong\Delta(2,3, 2n)$, the subgroups $K$ in Case~4 lift back to normal subgroups of level $4n$ in $\Gamma$, with $\Delta_1$ and $K_0$ lifting back to the principal congruence subgroups $\Gamma(2)$ and $\Gamma(4)$ of levels $2$ and $4$, and $\Delta_2$ lifting back to a conjugate of $\Gamma_0(4)$ (see~\cite[Ch.~6]{JS87}).

\medskip

\noindent{\bf Example 16.} We can take $K$ to be the image in $\Delta^*$ of the principle congruence subgroup $\Gamma(4n)$ of level $4n$ in $\Gamma$. We have
\[|\Gamma:\Gamma(m)|=\frac{m^3}{2}\prod_{p|m}\left(1-\frac{1}{p^2}\right)\]
for each $m>2$, with $p$ ranging over the distinct primes dividing $m$, so the resulting dessins ${\cal H}_i$ satisfy
\[|{\rm Aut}\,{\cal H}_i|=|\Delta_i:K|=\frac{1}{6}|\Delta^*:K|=\frac{1}{6}|\Gamma:\Gamma(4n)|
=4n^3\prod_{2\ne p|n}\left(1-\frac{1}{p^2}\right).\]
The Riemann surface $X$ underlying these dessins is the modular curve $X(4n)$ associated with $\Gamma(4n)$. By equation~(2) in Case~3, it has genus
\[g=1+\frac{(2n-3)}{4n}|\Delta_i:K|=1+(2n-3)n^2\prod_{2\ne p|n}\left(1-\frac{1}{p^2}\right).\]
More generally we could take $K$ to be the image in $\Delta^*$ of any normal subgroup $N$ of $\Gamma$ such that $\Gamma(4)\ge N\ge \Gamma(4n)$. These all correspond to normal subgroups of $\Gamma/\Gamma(4n)\cong PSL_2({\mathbb Z}_{4n})$; McQuillan~\cite{McQ} has classified the normal subgroups of $PSL_2({\mathbb Z}_m)$ for all $m$.

\medskip

Example~16 is not typical, in the sense that `most' normal subgroups $N$ of finite index in $\Gamma$ are non-congruence subgroups~\cite{Jon86}. If $N$ is any normal subgroup of finite index in $\Gamma$, and $l$ is its level, then $N\cap \Gamma(4)$ is a normal subgroup of finite index and level $4n={\rm lcm}\,(4,l)$ in $\Gamma$, so we can take $K$ to be its image in $\Delta^*$. If $\Gamma/N$ has a non-abelian composition factor not of type $L_2(p)$ for any prime $p$, then $N$ is a non-congruence subgroup, and hence so is $N\cap\Gamma(4)$, since the non-abelian composition factors of $PSL_2({\mathbb Z}_m)$ all have this type.

\medskip

\noindent{\bf Example 17.} A Hurwitz group $H$ is a non-trivial finite quotient of $\Delta(2,3,7)$ (for instance, see Examples~1, 2 and 3), and hence a quotient $\Gamma/N$ of $\Gamma$ where $N$ has level $l=7$. By letting $K$ be the image of $N\cap \Gamma(4)$ in $\Delta^*$, we obtain pairs of hypermaps ${\cal H}_i$ of types $(14, 14, 14)$ and $(7, 28, 28)$ on the same surface $X={\mathbb H}/K$, with ${\rm Aut}\,{\cal H}_i\cong V_4\times H$ or $C_4\times H$ for $i=1, 2$. For example, Ree's family of simple groups $Re\,(3^e)$ for odd $e>1$ are all Hurwitz groups~\cite{Jon94, Mal}, as is the Monster simple group~\cite{Wil}; these all correspond to non-congruence subgroups $N$ of $\Gamma$.

\section{Appendix: counting formulae}

Finding a regular dessin of type $(l,m,n)$ with automorphism group $G$ is equivalent to finding a smooth epimorphism from the triangle group $\Delta(l, m, n)$ onto $G$, and hence to finding a generating triple $x, y, z$ of orders $l, m$ and $n$ in $G$, with $xyz=1$. One can count such dessins by using character theory to count solutions of this equation in $G$.

More generally, if $X_1, \ldots, X_r$ are conjugacy classes
in a finite group $G$, then the number of $r$-tuples $(x_1,\ldots, x_r)\in X_1\times\cdots\times
X_r$ such that $x_1\ldots x_r=1$ in $G$ is given by Frobenius's formula
\[{|X_1|.\,\ldots\, .|X_r|\over |G|}\sum_{\chi}{\chi(x_1)\ldots\chi(x_r)\over\chi(1)^{r-2}},\]
where $x_i\in X_i$ and $\chi$ ranges over the irreducible complex characters of $G$ (see~\cite{Fro}, \cite{Jon95} or~\cite[Theorem~7.2.1]{Ser}). In the particular case $r=3$ the number of triples $(x, y, z)$ chosen from conjugacy classes $X, Y, Z$, with $xyz=1$, is 
\begin{equation}
{|X|.|Y|.|Z|\over |G|}\sum_{\chi}{\chi(x)\chi(y)\chi(z)\over\chi(1)}.
\end{equation}

In order to have a smooth homomorphism, we choose the conjugacy classes $X, Y$ and $Z$ to consist of elements of the same orders $l, m$ and $n$ as the canonical generators of $\Delta$. If one can show that some triple $(x, y, z)$ generates $G$ (for instance, by showing that no maximal subgroup contains $x, y$ and $z$), then we have an epimorphism $\Delta\to G$. 
A more sophisticated approach, due to P.~Hall~\cite{Hal}, uses M\"obius inversion in the subgroup lattice of $G$ to enumerate generating sets and hence epimorphisms; see~\cite{DJ, Jon95} for some applications.
 
For any groups $\Delta$ and $G$, two epimorphisms $\Delta\to G$ have the same kernel if and only if they differ by an automorphism of $G$; since ${\rm Aut}\,G$ acts fixed-point-freely on generating sets, it has orbits of length $|{\rm Aut}\,G|$ on these epimorphisms, so one can count kernels $K$ by dividing the number of epimorphisms by $|{\rm Aut}\,G|$. Applying this to $\Delta=\Delta(l,m,n)$ gives the number of regular dessins of type $(l,m,n)$ with automorphism group $G$. See~\cite{Jon94} for an example of this technique, applied to the Ree groups $Re(3^e)$ as Hurwitz groups.

\end{document}